\documentclass[11pt]{amsart}
\usepackage{amsmath,amsfonts,amssymb,amsthm,epsfig,color}
\newcommand{\SL}{\operatorname{SL}}
\newcommand{\GL}{\operatorname{GL}}

\newcommand{\PSL}{\operatorname{PSL}}
\newcommand{\bb}{\mathbb}

\newcommand{\Z}{\bb Z}
\newcommand{\R}{\bb R}
\newcommand{\N}{\bb N}

\newcommand{\A}{\bb A}

\newcommand{\B}{\mathcal{B}}

\newcommand{\bk}{\mathbf{k}}
\newcommand{\Q}{\mathbb Q}
\newcommand{\cO}{\mathcal{O}}
\newcommand{\bbG}{\mathbb G}

\newcommand{\fs}{\mathbb{F}_s}
\newcommand{\ft}{\mathbb{F}_{s}((X^{-1}))}
\newcommand{\fin}{\mathbb{F}_{s}[X]}

\newcommand{\ve}{\mathbf{x}}
\newcommand{\vc}{\mathbf{v}}
\newcommand{\q}{\mathbf{q}}
\newcommand{\pe}{\mathbf{p}}
\newcommand{\app}{\mathcal{W}(\psi)}

\newcommand{\Mat}{\operatorname{Mat}}
\newcommand{\csp}{\operatorname{Cusp}}

\newcommand{\hcal}{\mathcal{H}}
\newcommand{\kbar}{\bar{k}}
\newtheorem{Theorem}[equation]{Theorem}

\newtheorem{lemma}[equation]{Lemma}
\newtheorem*{lemma*}{Lemma}

\newtheorem*{theorem*}{Theorem}
\newtheorem{Def}[equation]{Definition}
\numberwithin{equation}{section}

 \email{jayadev.athreya@yale.edu}
 \email{A.Ghosh@uea.ac.uk}
\email{amri@imsc.res.in}
\begin{document}
\title{Ultrametric Logarithm Laws I.}
\subjclass{37A17, 37D40, 11K60}
 \keywords{
Shrinking target properties, logarithm laws, Diophantine approximation, local fields}
\thanks{J.\ S.\ A. gratefully acknowledges support from NSF  grant DMS
    0603636. A.\ G.\ gratefully acknowledges support from NSF grant DMS 0700128 and an EPSRC fellowship}
\maketitle
%% Enter the first author's name and address:
\centerline{\scshape J.\ S.\ Athreya}
\medskip
{\footnotesize
 %% please put the address of the first author
 \centerline{Department of Mathematics, Princeton University\footnote{Current Address: Mathematics Department, Yale University, PO Box 208283, New Haven, CT 06520-8283}}
   \centerline{ Fine Hall, Washington Road
Princeton NJ 08544-1000 USA}
} %% Do not forget to end the {\footnotesize by the sign }

\medskip

\centerline{\scshape A.\ Ghosh }
\medskip
{\footnotesize
 %% please put the address of the second author
 \centerline{Department of Mathematics, The University of Texas at Austin\footnote{Current Address: School of Mathematics, University of East Anglia, Norwich, NR4 7TJ UK}}
   \centerline{
1 University Station, C1200
Austin, Texas 78712}
} %

\medskip

\centerline{\scshape A.\ Prasad }
\medskip
{\footnotesize
 %% please put the address of the second author
 \centerline{Department of Mathematics, Institute of Mathematical Sciences}
   \centerline{
Taramani, Chennai 600 113, India.}
} %

\bigskip

%% The name of the associate editor will be entered by an editorial staff
% \centerline{(Communicated by the associate editor name)}

\begin{abstract} 
We announce ultrametric analogues of the results of Kleinbock-Margulis for shrinking target properties of semisimple group actions on symmetric spaces. The main applications are $S$-arithmetic Diophantine approximation results 
and logarithm laws for buildings, generalizing the work of Hersonsky-Paulin on trees.
\end{abstract}

\section{Introduction}\label{intro} 
The notion of shrinking targets for dynamical systems is a much studied and extremely useful concept (\cite{Athreyasurvey}, \cite{KC} and the references therein) especially for geometric and number theoretic applications. Shrinking target properties for group actions on homogeneous spaces were studied in an important paper \cite{KMinv} of D.\ Kleinbock and G.\ A.\ Margulis. We first describe their results. 
Let $G$ denote a connected semisimple Lie group without compact factors, and $\Gamma$ denote a non-uniform lattice in $G$. Let $K$ be a maximal compact subgroup of $G$, and let $Y = K \backslash G/\Gamma$ denote the associated non-compact irreducible locally symmetric space of finite volume. In \cite{KMinv},  D.\ Kleinbock and G.\ A.\ Margulis studied the phenomenon of shrinking target properties for the geodesic flow, generalizing an earlier work of D.\ Sullivan \cite{Sullivan} and established the following theorem, commonly called a logarithm law.\\

 For $x \in Y$, let $T^{1}_x(Y)$ denote the unit tangent space at $x$. Let $\nu$ denote the Haar measure on $T^{1}_x(Y)$. For $\theta \in T^{1}_x(Y)$, and $t \in \R$ let $g_t(x, \theta)$ denote the image of $(x, \theta)$ under geodesic flow for time $t$. Let $d_Y$ denote a metric on $Y$, obtained from a bi-$K$ invariant metric $d$ on $G$.

\begin{Theorem}\label{loglawzero}(Kleinbock-Margulis~\cite{KMinv}) There exists a $k = k(Y)>0$ such that the following holds: for all $x, y\in Y$, and almost every $\theta \in T^1_x(Y)$, $$\limsup_{t \rightarrow \infty} \frac{d_Y(g_t(x, \theta), y)}{\log t} = 1/k.$$\end{Theorem}

The Theorem thus studies the statistical properties of geodesic excursions into smaller and smaller cuspidal neighborhoods (i.e. the complements of these neighborhoods are larger and larger compacta) of $Y$. In fact, this phenomenon seems to be prevalent in many dynamical systems, and a very general result to the effect was obtained in \cite{KMinv}, which can then be used to obtain Theorem \ref{loglawzero}. To describe this result, we need some more definitons,  all taken from loc. cit. Let $(X,\mu)$ be a probability space, $\B$ be a family of measurable subsets of $X$ and let $F = \{f_t\}$ denote a sequence of $\mu$-preserving transformations of $X$. Then, 

\begin{Def}
$\B$ is called Borel-Cantelli for $F$ if for every sequence $\{A_t~|~t \in \N\}$ of sets from $\B$, the following holds:
\begin{displaymath}
\mu(\{x \in X~|~f_{t}(x) \in A_t~\text{for infinitely many}~t \in \N\})
\end{displaymath}
\begin{displaymath}
= \left \{ \begin{array}{ccc}
0 & \textrm{if $\sum_{t=1}^{\infty}\mu(A_t)$ converges}, \\\\
1 & \textrm{if $\sum_{t=1}^{\infty}\mu(A_t)$ diverges}.
\end{array}\right.
\end{displaymath}

\end{Def}

For a function $\delta$ on $X$, denote by $\Phi_{\delta}$, the tail distribution function, defined by:
\begin{equation}
\Phi_{\delta}(z) = \mu(\{x~|~\delta(x) \geq z\}).
\end{equation}

\noindent And for $\kappa > 0$, we say that $\delta$ is $\kappa-DL$ (an abbreviation for $\kappa-$distance like) if it is uniformly continuous and 
\begin{equation}\label{kDL}
\exists~C_1, C_2 > 0,~\text{such that}~C_1e^{-\kappa z}\leq \Phi_{\delta}(z)\leq C_2e^{-\kappa z}~~\forall ~z \in \R,
\end{equation}

\noindent and $DL$ if there exists $\kappa > 0$ such that (\ref{kDL}) holds. The following is then Theorem $1.8$ in \cite{KMinv}.

\begin{Theorem}\label{km1.8}
Let $G$, $\Gamma$ and $X$ be as above, $F = \{f_t~|~t \in \N\}$ be a sequence of elements of $G$ satisfying
\begin{equation}\label{ED}
\sup_{t \in \N} \sum_{s=1}^{\infty} e^{-\beta d(f_sf^{-1}_t,e)} < \infty~\forall~\beta > 0,
\end{equation}

\noindent and let $\delta$ be a $DL$ function on $X$. Then 
\begin{equation}\B(\delta)\overset{def}= \left\{ \{ x \in X~|~\delta(x) \geq r\}~|~r \in \R \right\}\end{equation} 
\noindent is Borel-Cantelli for $F$.

\end{Theorem}

\noindent Sequences which satisfy (\ref{ED}) above are referred to as ``exponentially divergent" (abbreviated $ED$) in \cite{KMinv}. To derive Theorem \ref{loglawzero} from Theorem \ref{km1.8}, the authors use a philosophy of F.Mautner \cite{Mautner} to realize the geodesic flow on the unit tangent bundle of $Y$ as a one-parameter flow on $G/\Gamma$. One then needs to check the $ED$ condition for this flow, as well as the $DL$ condition for the metric on $G/\Gamma$ described above. In this paper and its sequel \cite{agp2}, we will obtain $S$-arithmetic analogues of Theorem \ref{km1.8}.\\ 

To motivate our results, we start with an analogue of Theorem \ref{loglawzero}.  Let $\fs$ denote the finite field of $s = p^{\vartheta}$ elements, $\bk$ denote the global function field of rational functions with coefficients in $\fs$ and $k$ denote the completion of $\bk$ at the infinite place, identified naturally with the field of Laurent series $\ft$ with coefficients in $\fs$. On $\ft$, we will assume the usual norm and ultrametric, (cf.\cite{Weil}) for details. We denote this norm $\|~\|$ and when used in the context of $k^r$, it will refer to the supremum norm. Let $\bbG$ denote a simple, isotropic, linear algebraic group defined and split over $\fs$, $\Gamma$ a non-uniform lattice in $\bbG(k)$\footnote{It is well-known $G$ contains many such lattices.}, and $K$ denote a parahoric subgroup of $\bbG(k)$. We fix a bi $K$-invariant metric on $\bbG(k)$, and hence also on $\bbG(k)/\Gamma$ and call this metric $d(~,~)$. When the field in question is evident, we will sometimes refer to $\bbG(k)$ simply as $G$.\\  

A natural geometric object on which $G$ acts is the Bruhat-Tits building $Y$ of $G$. This is a Euclidean building which comes with a natural metric, equipped with which it becomes a CAT-$0$ space. The stabilizers of vertices in $Y$ are the maximal parahoric subgroups $K$ of $G$. Therefore the quotient $Y/\Gamma$ can be naturally identified with $K\backslash G/\Gamma$.  Let $\partial Y$ denote the geodesic boundary of $Y$ and $\pi: Y \rightarrow Y/\Gamma$ be the natural projection. Let $\{ g_t (x, \theta)\}_{t \geq 0}$ denote the geodesic starting at  the vertex $x \in Y$ in direction $\theta \in \partial Y$. There is a natural measure class on $\partial Y$. A function field analogue of Theorem \ref{loglawzero} would therefore be:

\begin{Theorem}\label{loglaw} There is a $l = l(Y/\Gamma)$ such that  for any $x \in Y$, $y \in Y/\Gamma$ and almost all $\theta \in \partial Y$, $$\limsup_{t \rightarrow \infty} \frac{ \log d_{Y/\Gamma} (\pi(g_t(x, \theta)), y)}{ \log t} = 1/l,$$\noindent where $d_Y$ is the metric on the quotient $Y/\Gamma$.\end{Theorem}

\noindent\textbf{Remark:} This result was obtained by S.\ Hersonsky and F.\ Paulin in~\cite{Her-Pau} in the case of quotients of a regular tree $Y$ by any non-uniform lattice in $Aut(Y)$.  
Their approach is more geometric, since $Aut(Y)$ is a locally compact group but does not have any obvious algebraic structure, and in some sense is in the spirit of Sullivan \cite{Sullivan}. A natural
question would be to ask if an analogue of Theorem \ref{km1.8} holds in a more general algebraic setting. We answer this in the affirmative.\\

Let $\bk$ be a global field, $S$ a finite set of places of $\bk$ (containing the infinite ones, in case $\bk$ is a number field), and for each $s \in S$, let $k_s$ denote the completion of $\bk$ at the place $s$,  $\bbG$ be a connected, simple, algebraic $\bk$-group without anisotropic factors, set $G_s = \bbG(k_s)$ and $G _S= \prod_{s \in S}G_s$. Let $Y_s$ denote the Bruhat-Tits building of $G_s$ (or the symmetric space as the case may be) and $Y_S = \prod_{s \in S}Y_s$. Let $\Gamma$ be a non-uniform lattice in $G_S$ and set $X_S = G_S/\Gamma$. In \cite{agp2}, we prove:

\begin{Theorem}\label{agp1}
With notation as above, let $F$ be an exponentially divergent sequence for $G_S$ and $\delta$ be a $DL$ function on $X_S$.
 Then $\B(\delta)$ is Borel-Cantelli for $F$. 
\end{Theorem}

\noindent\textbf{Remark 1.} The reader will notice that Theorem \ref{km1.8} is stated in somewhat greater generality, i.e. in the context of semisimple groups. However, the principal ideas involved are already present in the case of simple groups. In \cite{agp2} where the proof of the above theorem will be presented in detail, we will point out the extra details required to extend Theorem \ref{agp1} to semisimple groups.\\
\noindent\textbf{Remark 2.} As a corollary, we will obtain generalizations of Theorems \ref{loglawzero} and \ref{loglaw}  for $S$-isotropic symmetric spaces.\\

We now turn to some examples : Let $\cO_S$ be the ring of $S$-integers of $\bk$, i.e. 
\begin{equation}
\cO_S = \{x \in \bk~|~|x|_s \leq 1~\text{for every finite place}~s \notin S\}.
\end{equation}

Then, by theorems of Borel-Harish Chandra and Behr-Harder in the function field case (cf.\cite{Marg}),  $\Gamma \overset{def}= \bbG(\cO_S)$ is a lattice\footnote{This, and commensurable lattices are called arithmetic.} in $G_S$. Let $K_s$ denote a maximal parahoric subgroup (maximal compact if $s$ is an infinite place) of $G_s$, and $d_s$ denote the $K_s$-bi invariant metric on $G_s$. This gives rise to a product metric $d$ on $X_S$. In \cite{agp2}, using reduction theory we show that $d$ is a $DL$ function on $X_S$ for \emph{any} non-uniform lattice $\Gamma$, thereby establishing Theorem \ref{loglaw} by an application of Theorem \ref{agp1}.\\

 For an application to number theory, take $G = \SL(n,\R), \Gamma = \SL(n,\Z)$, so the space 
 \begin{equation}
 X = \SL(n,\R)/\SL(n,\Z)
 \end{equation}
\noindent  can be identified with the space of unimodular lattices in $\R^n$. In \cite{KMinv}, Theorem \ref{km1.8} and the function
\begin{equation}
\delta(\Lambda)=\inf_{\ve  \in \Lambda \backslash \{0\}}\| \ve \|
\end{equation}
 which turns out to be $DL$, was used to prove Khintchine's theorem, a cornerstone of metric Diophantine approximation, and stronger variants of the it. Now let $\bk = \Q$ and  following \cite{KleTom}, define

\begin{equation}
\GL^{1}(n, \Q_S)\overset{def}= \left\{g = (g_s)_{s \in S} \in \GL(n,\Q_S)~|~\prod_{s \in S}\det(g_s)_{s \in S} = 1\right\}
\end{equation} 

Then $\GL^{1}(n,\Q_S)/\GL(n, \Z_S)$ is the space of unimodular lattices in $\Q_S$. In \cite{KleTom}, a dynamical interpretation of Diophantine approximation, originally introduced by D.\ Kleinbock and G. \ Margulis in \cite{KM2} was developed in the $S$-arithmetic setting by D.\ Kleinbock and G.\ Tomanov to prove $S$-arithmetic analogues of Mahler's conjectures, see also \cite{G} for function field analogues. In \cite{agp2}, we use Theorem \ref{agp1} in this context with specific, analogous choices of $DL$ functions to establish Khintchine's theorems  and multiplicative versions in the number and function field cases. We illustrate with an example. Let $l$ denote the cardinality of $S$, and assume, for simplicity that $S$ contains the infinite valuation\footnote{This is just a convenience because the definitions are slightly more involved cf.\cite{KleTom}. This restriction will be removed in \cite{agp2}}. Let $\psi : \N \to \R_{+}$ be a non-increasing function, and let $\Mat(\Q_S)$ denote the set of $m \times n$-matrices with entries in $\Q_{S}$. Let $\app$ denote the set of $A \in \Mat(\Q_S)$ for which there exist infinitely many $\q \in \Z^n$ such that
\begin{equation}\label{psiapp}
(\|\pe + A\q\|^{l})^{m}\leq \psi(\|\q\|^{n}_{\infty})~\text{for some} ~\pe \in \Z^{m}.
\end{equation}
 Here the norms $\|~\|$ and $\|~\|_{\infty}$ are defined as follows:\\
\noindent For $\ve = (x^{(v)}_1,x^{(v)}_2,\dots,x^{(v)}_m) \in \Q^{m}_v$, we define 
\begin{equation}\|\ve\|_v = \max_{i}|x^{(v)}_{i}|_{v}
\end{equation} 
\noindent and finally, for $\ve \in \Q^{m}_{S}$
\begin{equation}
\|\ve\| = \max_{v}\|\ve\|_v.
\end{equation}

For a very nice motivation of the theory of $S$-arithmetic Diophantine approximation, and especially the correct normalization, we refer the reader to Section $10$ in \cite{KleTom}. The $S$-arithmetic analogue of the Khintchine-Groshev theorem would then be:
 \begin{Theorem}\label{SKG}
 Let $\app$ denote the set of $\psi$-approximable matrices as above. Then, 
 \begin{enumerate}
 \item Almost every $A$ belongs to $\app$ if $\int\psi(x)~dx$ diverges.\\
 \item Almost no $A$ belongs to $\app$ if $\int\psi(x)~dx$ converges.
 \end{enumerate}
 \end{Theorem}  
 This theorem will be proved in \cite{agp2} by adapting the method of Kleinbock-Margulis to this set-up. In addition, we will prove variations and strengthenings of the above theorem, including:
 \begin{enumerate}
 \item The more general multiplicative case, where the norm $\|~\|$ is replaced by a more general function. This will prove $S$-arithmetic analogues of a conjecture of Skriganov \cite{Sk}.\\
 \item Function field analogues of Theorem \ref{SKG}. This is elaborated upon in the next section. 
 \end{enumerate}

\noindent This paper is intended as an announcement of our results above.  In particular, we will not present complete proofs, but rather try to convey the main ideas involved in the proof. To do this, we will focus on the case of a single local field of positive characteristic. Namely,  we will outline the proofs of  Theorems \ref{loglaw} and \ref{SKG} in the case where $k$ is a local field of positive characteristic. We would like to stress that our strategy is substantially similar to that of Kleinbock-Margulis. The primary content of our work is in generalizing their methods to much wider settings. However, in \cite{agp2}, we will also provide a direct, ``non-spectral" proof of Theorem \ref{loglaw} in the case where $\bbG$ is a rank $1$ group, which has been explained to us by S. Mozes. In this case, we use a symbolic description of the geodesic flow on trees, developed in \cite{Mozes}. We are indebted to him for sharing his insight.

\subsection*{Acknowledgements.}  A.\ G. announced these results at an AMS regional conference in Chicago, Ill. and would like to thank Marian Gidea for his hospitality. J.\ A. and A.\ G. thank D.\ Kleinbock and the participants of the Workshop on Shrinking Target Properties, Feb. 2008 for helpful conversations. We also thank  the Clay Mathematical Institute for sponsoring the workshop. We thank the referee for several helpful comments which have improved the exposition of the paper.

\section{Borel-Cantelli lemmata and applications.}

Theorem \ref{agp1}, like Theorem \ref{km1.8} are probabilistic in nature-in fact they are strongly reminiscent of $0-1$ laws in probability theory, especially the elementary Borel-Cantelli lemma, which we now recall: 
\begin{lemma}\label{BC}(Borel-Cantelli) Let $\{X_n\}_{n=0}^{\infty}$ be a sequence of $0-1$ random variables, with $P(X_n = 1) =: p_n$. Then 
\begin{enumerate}
\item If $\sum_{n=0}^{\infty} p_n < \infty$, then $P(\sum_{n=0}^{\infty} X_n = \infty) =0$
\item If the $X_n$'s are pairwise independent ,i.e. 
\begin{equation}\nonumber
p_{nm} := P(X_n X_m =1) = p_n p_m~\forall~m, n, 
\end{equation}
\noindent and $\sum_{n=0}^{\infty} p_n = \infty$,  then $P(\sum_{n=0}^{\infty} X_n = \infty) =1$.
\end{enumerate}
\end{lemma}

The first part of the above Lemma easily allows one to derive the convergence halves of the various $0-1$ laws we have in mind. Unfortunately, it is typically hopeless to expect that the dynamical random variables (or events) we are interested in are independent, i.e. satisfy $(2)$ in the above theorem. In the context of logarithm laws, for instance, one is only able to show that geodesic excursions to shrinking cusp neighborhoods are relatively (sometimes referred to as ``quasi") independent events. That this suffices for applications is due to the following strengthening of the Borel-Cantelli lemma, abstracted from the works of W. Schmidt, by V. Sprindzhuk \cite{Sprindzhuk}. We first set up some of the notation, taken from \cite{KMinv}. For a function $f$ on a probability space $(X,\mu)$, we will denote
\begin{equation}\nonumber
\mu(f) \overset{def}= \int_{X} f~d\mu.
\end{equation}
\noindent and for $N \in \N \cup \{\infty\}$, a family of functions $\hcal = \{h_t~|~t \in \N\}$ on $X$,
\begin{equation}\label{sumavg}
S_{\hcal,N} \overset{def}= \sum_{t=1}^{N}h_t(x)~\text{and}~E_{\hcal,N}\overset{def}= \sum_{t=1}^{N}\mu(h_t).
\end{equation}

\begin{lemma}\label{Sprindzhuk}
Let $(X,\mu)$ denote a probability space, and let $\hcal = \{h_t~|~t \in \N \}$ denote a sequence of functions on $X$ which satisfy:
\begin{equation}\label{Sprindzhuk1}
\mu(h_t) \leq 1~\text{for every}~t \in \N.
\end{equation}
\noindent Assume also that there exists $C > 0$ such that
\begin{equation}\label{Sprindzhuk2}
\sum_{s,t = M}^{N}\left(\mu(h_sh_t)-\mu(h_s)\mu(h_t)\right) \leq C.\sum_{t=M}^{N}\mu(h_t)~\text{for every}~N > M \geq 1. 
\end{equation}
Then for every $\epsilon > 0$,
\begin{equation}\label{Sprindzhuk3}
S_{\hcal,N} = E_{\hcal,N}+O\left(E^{1/2}_{\hcal,N}\log^{3/2+\epsilon}E_{\hcal,N} \right)
\end{equation}
\noindent for $\mu$ almost every $x \in X$. In particular, for $\mu$ almost every $x$, 
\begin{equation}
\lim_{N \to \infty}\frac{S_{\hcal,N}(x)}{E_{\hcal,N}} = 1,
\end{equation}
\noindent whenever $\sum_{t=1}^{\infty}\mu(h_t)$ diverges. 
\end{lemma}

So, in order to use Sprindzhuk's lemma to prove Theorem \ref{agp1}, we need to show that the ``quasi-independence" condition (\ref{Sprindzhuk2}) holds in the context of our dynamical systems. We start with the Khintchine-Groshev theorem and its generalizations, focussing on the function field case.  Our plan is to first describe a dynamical system to which the above Lemma can be applied to derive the desired number theoretic results, and then to use quantitative mixing bounds for the dynamical system to ensure that condition (\ref{Sprindzhuk2}) is satisfied.\\
 
Accordingly, let us take $G = \SL(n,\ft)$ and $\Gamma = \SL(n,\fin)$. Then $\Gamma$ is a non-compact lattice in $G$. Let $X = G/\Gamma$ and $\mu$ denote the normalized, projected Haar measure on $X$. Diophantine approximation in the function field setting is naturally analogous\footnote{This is true to some extent. There are important and much studied distinctions, especially in the approximation of algebraic elements. We will not address this.} to the real case, namely the role of $\R$ is played by $\ft$, while that of $\Z$ is played by the polynomial ring $\fin$, and there is an analogous continued fraction decomposition for every element of $\ft$ (see, for example, ~\cite{Paulin}). As would be expected, it is possible to read off Diophantine properties of Laurent series from their continued fraction expansions, and in fact one can provide a proof of Khintchine's theorem in one dimension using a careful study of these continued fractions.\\

More generally, the set $\app$ of $\psi$-approximable $(m \times n)$ matrices with entries in $\ft$ is naturally defined to be those $A \in \Mat_{m,n}(\ft)$ for which there exist infinitely many $\q \in \fin^{n}$ such that
\begin{equation}\label{Diophmat}
\|A\q+\pe\|^{m}<\psi(\|\q\|^{n})~\text{for some}~\pe \in \fin^{m}.
\end{equation}

\noindent And the function field analogue of the Khintchine-Groshev Theorem would precisely be Theorem \ref{SKG} with the above definition of $\app$.\\

An ingenious scheme, due to Dani \cite{Dani} in a special case, and developed in full generality by Kleinbock-Margulis \cite{KMinv} translates the above problem into one of shrinking target properties for certain homogeneous flows. Let us briefly describe this correspondence. The group $\SL_{m+n}(\ft)$ acts transitively on the space of unimodular (i.e. co-volume $1$) lattices $\Omega_{m+n}$ of $\ft^{m+n}$ and the stabilizer of $\fin^{m+n}$ is $\SL_{m+n}(\fin)$.\\ 

The space $\SL_{m+n}(\ft)/\SL_{m+n}(\fin)$ can thus be naturally identified with $\Omega_{m+n}$. Given $A \in \Mat_{m,n}(\ft)$, we associate to it the following lattice:
\begin{equation}\nonumber
A \rightsquigarrow \Lambda_{A}\overset{def}=\begin{pmatrix}I_m & A \\ 0 & I_n\end{pmatrix}\fin^{m+n}
\end{equation}
\noindent where $I_i$ denotes the square identity matrix of dimension $i$. Good rational approximations to $A$ can be shown to correspond to small vectors in the lattice $\Lambda_{A}$ corresponding to $A$. Let 
 $\pi$ denote the uniformizer of $k$ and  set
\begin{equation}
g_{t} \overset{def}= diag(\underbrace{\pi^{nt},\dots, \pi^{nt}}_{m~\text{times}}\underbrace{\pi^{-mt}, \dots, \pi^{-mt}}_{n~\text{times}}).
\end{equation}
\noindent and let $\delta$ denote the following function which measures small vectors in lattices:
\begin{equation}\label{defdelta}
\delta = \min_{\vc \in \Lambda \backslash\{0\}}\| \vc \|.
\end{equation}

By the above identification of \begin{equation}\nonumber\Omega_{m+n}~\text{with}~\SL(m+n,\ft)/\SL(m+n,\fin)\end{equation} and Mahler's compactness criterion, those $t$ for which $\delta(g_t\Lambda_{A})$ is small correspond to $t$ for which the $g_t$ trajectory of $\Lambda_A$ has cusp excursions. To define these excursions, we set
\begin{equation}
\csp(t,\delta) \overset{def}=\{\Lambda \in \Omega_{m+n}~|~\delta(\Lambda)\leq |\pi|^{-t} \}.
\end{equation}

\noindent The following lemma, proven\footnote{In the real case, but the function field proof goes through with obvious modifications. An $S$-arithmetic version will be provided in \cite{agp2}.} in \cite{KMinv} then establishes the precise connection between Diophantine properties and cusp excursions.
\begin{lemma}\label{cov} There exists an explicit function $r(t)$ which depends only on $\psi, m$ and $n$ and satisfies
\begin{equation}\nonumber
\int \psi(x)~dx < \infty \iff \int |\pi|^{-(m+n)r(t)}~dt < \infty.
\end{equation}
Moreover, $A \in \app$ if and only if
\begin{equation}\label{cuspexc} \exists~\text{infinitely many}~t \in \N~\text{such that}~g_{t}\Lambda_{A} \in \csp(r(t),\delta).\end{equation}
\end{lemma}

So in order to prove Theorem \ref{SKG} for function fields, it clearly suffices to prove that:
\begin{Theorem}\label{KGlattices} With notation as above,\\
\begin{enumerate}
\item Almost every $\Lambda_A$ satisfies (\ref{cuspexc}) if $\int |\pi|^{-(m+n)r(t)}~dt $ diverges.\\
\item Almost no $\Lambda_A$ satisfies (\ref{cuspexc}) if $\int |\pi|^{-(m+n)r(t)}~dt$ converges.
\end{enumerate}
\end{Theorem}

Similarly, in order to prove logarithm laws for $K\backslash \bbG(k)/\Gamma$,  where $\bbG$ is a simple group as before, we will use the philosophy of F. Mautner \cite{Mautner} to realize the geodesic flow on the unit tangent bundle of $K\backslash \bbG(k)/\Gamma$ as a one-parameter flow $g_t$ on $\bbG(k)/\Gamma$. Fix $x_0 \in \bbG(k)/\Gamma$,  and set 
\begin{equation}
\csp(t,d)=\{y \in \bbG(k)/\Gamma~|~d(x_0,y) > t\}.
\end{equation}
The object of study then becomes the excursions of $g_t$ orbits into $\csp(t,d)$.\\

We now tie these themes up with Sprindzhuk's lemma, focussing for convenience on Theorem \ref{KGlattices}. Accordingly, let $h_t$ denote the characteristic function of $\csp(r(t),\delta)$, and let 
\begin{equation}
\hcal^{\mathcal{G}} = \{g^{-1}_{t}h_t~|~t \in \N\}.
\end{equation} 
Then to derive Theorem \ref{KGlattices} from Sprindzhuk's Lemma, we need to ensure that
\begin{equation}
\sum_{s,t=M}^{N}\left((g^{-1}_{s}h_s,g^{-1}_{t}h_t)-\mu(h_s)\mu(h_t)\right)
\end{equation}

\noindent is small, in other words we want to show that the excursions of $g_t$ orbits into $\csp(r(t),\delta)$ are quasi-independent. The key to this is the spectral gap of the $G$-action on $G/\Gamma$. 

\section{Spectral Gap.}
We retain the notation of the previous section, i.e. $G$ is the group of $k$-points of a simple $\bbG$ as before, and $\Gamma$ is a non-uniform lattice in $G$. A natural tool to study the ergodic properties of the action of $G$ on $G/\Gamma$ is the spectral properties of the action of $G$ on $L^{2}(G/\Gamma)$. Let $L^{2}_{0}(G/\Gamma)$ denote the subspace of $G/\Gamma$ orthogonal to constants. Let $K$ denote a maximal compact open subgroup of $G$. We call $\phi \in L^{2}_{0}(G/\Gamma)$, $K$-finite if its $K$-span is finite dimensional. On $G$ there is  an important function which controls the rate of decay of matrix coefficients, the Harish-Chandra function $\Xi$.  Let  $G = KAN$ denote the Iwasawa decomposition of $G$. The Harish-Chandra function of $G$ is then defined by:
\begin{equation}\label{HC}
\Xi(g) = \int_{K} \delta^{-1/2}(gk)~dk
\end{equation}

\noindent where $\delta$ is the left modular function of $AN$ defined by:
\begin{equation}\label{delta}
d\mu_{G}=d\mu_{K}\delta(an)d\mu_{AN}
\end{equation}

\noindent where $\mu_{*}$ denotes Haar measure on the locally compact group $*$. The following estimate for decay of matrix coefficients is then known to hold by work of several authors. We refer to Section 5.1 in \cite{GN} as a convenient reference.
\begin{Theorem}\label{matdec1}
Let $G, K$ and $\Gamma$ be as above. Then there exist positive constants $C,\chi$ such that for every $K$-finite $\phi, \psi \in L^{2}_{0}(G/\Gamma)$ and any $g \in G$, 
\begin{equation}
|<g\phi,\psi>| \leq C\|\phi\|\|\psi\|\Xi(g)^{1/\chi}.
\end{equation}
\end{Theorem}

It turns out that the above estimate, or more precisely a version of this estimate which caters to smooth functions (the $L^2$-norm is then replaced by an appropriate Sobolev norm) is enough to ensure the quasi-independence condition in Sprindzhuk's lemma. The passage from $L^2$ to smooth functions is quite standard and will be elaborated upon in the $S$-arithmetic setting in \cite{agp2}. It remains to show that the characteristic functions of $\csp(r(t),\delta)$ (resp. $\csp(t, d)$) can be approximated by appropriate smooth functions, which is precisely showing that the functions in question are $DL$. This ``smoothing of cusp neighborhoods" is carried out using reduction theory, which we elaborate upon in the next section. \\

We return briefly to the proof of Theorem \ref{matdec1}. It turns out that these bounds are closely related to the notion of spectral gap for the $G$ action.
Recall that the $G$ action on a probability space $X$ has spectral gap if the regular representation $\rho_0$ of $G$ on the space $L^{2}_{0}(X)$, the subspace of $L^{2}(X)$ orthogonal to constants, is isolated in the Fell topology from the identity (or trivial) representation. If $G$ is an almost direct product of groups $G_i$, then the $G$ action on $X$ has strong spectral gap if the restriction of $\rho_0$ to any factor $G_i$ is isolated from the trivial representation.\\ 

Establishing strong spectral gap turns out to be quite difficult in general. In \cite{KMinv}, the authors proved that if $G$ is a connected semisimple, center-free Lie group without compact factors, $\Gamma$ is an irreducible non-uniform lattice in $G$, then the $G$ action on $G/\Gamma$ has strong spectral gap. The analogous question for uniform lattices does not seem to be known, however see the recent work \cite{KS}. Of course, if all the factors $G_i$ have Kazhdan's property T, strong spectral gap is immediate.\\

 However, since rank $1$-groups defined over non-Archimedean local fields act on their Bruhat-Tits trees without fixed points, they cannot have property $T$. A weaker and very useful property is property $\tau$ as defined by Lubotzky and Zimmer \cite{Lu-Zim}, which means that $\rho_0$ is isolated in the Automorphic Spectrum of $G$. In \cite{agp2}, we will gather various tools from representation theory namely the restriction technique of Burger-Sarnak \cite{BS}, as extended to finite places and function fields by Clozel-Ulmo (cf. \cite{Clo-Ull} and \cite{Clo}) and property $\tau$ for congruence subgroups of $\SL_2$ as established by Selberg \cite{Selberg}, Gelbart-Jacquet \cite{Gel-Jac}, Clozel \cite{Clotau} and Drinfeld \cite{Drinfeld} in various contexts, to record a proof of (the $S$-arithmetic generalization) of Theorem \ref{matdec1}.\\

We note that in many cases  uniform, optimal estimates for decay of matrix coefficients are known due to H.Oh \cite{Oh} and these have several important applications. While the use of spectral gap is a powerful tool, we remark that very little is known in the context of the automorphism group of a tree or more generally, a building and this poses a potential hurdle to obtain logarithm laws in these settings using spectral methods.\\

The reader will notice that while the existence of the constant $\chi$ in Theorem \ref{matdec1} is crucial as well as sufficient for the purposes of this paper, we do not say anything regarding bounds, i.e. the so-called extent of temperedness of the representation. These bounds are related to bounds towards the Generalized Ramanujan Conjecture and play a crucial role in many Diophantine problems. In \cite{GG}, the second named author and A.\ Gorodnik study the connection of Diophantine approximation on symmetric spaces and bounds on temperedness in greater detail. This problem involves studying shrinking target properties for targets which are balls around a fixed point in the symmetric space, and needs new techniques.
 
\section{Ad\`{e}lic Reduction Theory.}
We now present the tools required to obtain smoothing of the cusp neighborhoods, i.e. to show that the characteristic functions of super-level sets of the various cusp neighborhoods are $DL$. The main tool is ad\`{e}lic reduction theory, as developed in \cite{Pra01} and \cite{Pra03} from where we also borrow notation. Let $\bk, \bbG$ etc. be as before and let $\A$ denote the ad\`{e}le ring of $\bk$. Let $B$ be a fixed Borel subgroup of $\bbG(k)$ containing a maximal split torus $T$ of $G$. Let $X_{*}(T)$ denote the lattice of cocharacters of $T$. Let $X_{*}(T)^{++}$ denote the subset of $X_{*}(T)$ consisting of dominant cocharacters.  
Let $W = N_{G}T(\fs)/T(\fs)$ denote the Weyl group of $G$. Define the map $\phi_{\pi} : W \times X_{*}(T) \to G(\kbar)$ by
\begin{equation}\label{phidef2}
\phi_{\pi} (w,\mu) = w\mu(\pi).
\end{equation}
Let $\rho$ denote the sum of all roots that are positive with respect to $B$.  For $\tau > 0$, we denote by $X_{\tau}$ the subset of $\bbG(\A)/\bbG(k)$ consisting of the union of $\bbG(k)\phi_{t^{-1}}(e^{\lambda})K$ where $\lambda$ ranges over the dominant cocharacters of $T$ for which $\langle\rho,\lambda\rangle \geq \tau$. Then, the following estimate, which follows from the reduction theory developed in \cite{Pra01} and \cite{Pra03} allows us to show the requisite $DL$ property for $d(x_0,~)$.

\begin{Theorem}\label{adelicrec}
Let $r$ denote the $k$-rank of $\bbG(k)$. Then, there exist positive real numbers $C_1$ and $C_2$ such that for every $T >  0$,
\begin{equation}
C_1 < \frac{\mu(X_T)}{\sum_{l \geq T}s^{-l}l^{r-1}} < C_2.
\end{equation}
\end{Theorem}

We now turn to Khintchine's theorem , where in light of the previous discussion, we have to develop the reduction theory in a slightly different setting-namely on the space $\SL(d,\ft)/\SL(d,\fin)$ and for the function $\delta$. The key to this, as shown in \cite{KMinv}, lies in the so-called Siegel integral formula \cite{Siegel}, specifically to a multi-dimensional version of it.  We will again proceed in the ad\`{e}lic setting. Set:
\begin{equation}\label{adelic lattice} 
 \GL^{1}(d,\A)\overset{def}=\{g \in \GL(d,\A)~|~\prod_{v \in V(\bk)}|\det g|_{v} = 1\}
 \end{equation}
 
 \noindent where $V(\bk)$ denotes the set of inequivalent places of the global field $\bk$. Let $f \in L^{1}(\A^{d})$, $g \in \GL^{1}(d,\A)$ and define 
 \begin{equation}
 \tilde{f}\overset{def}= \sum_{\ve \in g\Q^d}f(\ve)
 \end{equation}
 
 \noindent Then the ad\`{e}lic version of Siegel's integral formula, as proved by Weil in \cite{Weil} states that:
 \begin{Theorem}\label{adelicSiegel}
 \begin{equation}\nonumber
 \int_{\A^d}f~d\mu_{\A} = C(d)\int_{\GL^{1}(d,\A)/\GL(d,\Q)} \tilde{f}~d\tilde{\mu}.
 \end{equation}
 \end{Theorem}
 \noindent where $C(d)$ is a constant depending on the field in question (in this case $\Q$). We remark that the above Theorem was proved by Weil for arbitrary number fields as well as function fields. A multidimensional generalization of this theorem is developed in \cite{agp2} and is shown to imply a local version of this formula which can then be used to show that the function $\delta$ is $DL$.

\end{document}